\documentclass[a4paper,english,oneside]{article}
\usepackage[T1]{fontenc}
\usepackage[utf8]{inputenc} 
\usepackage{authblk}

\usepackage{geometry}
\usepackage{amsmath}
\usepackage{amsfonts}
\usepackage{amssymb}
\usepackage{amsthm}
\usepackage[english]{babel}
\usepackage{color}
\usepackage{lmodern}
\usepackage{pstricks-add}
\usepackage{hyperref}
\hypersetup{hidelinks}

\usepackage{mathscinet}
\usepackage{enumitem}
\usepackage[cal=boondoxo,bb=ams]{mathalfa}
\usepackage{microtype}

\usepackage{mathtools}
\usepackage{esint}

\theoremstyle{plain}
\newtheorem{theorem}{Theorem}[section]

\theoremstyle{definition}
\newtheorem{definition}[theorem]{Definition}

\theoremstyle{remark}
\newtheorem{remark}{Remark}

\begin{document}

\title{On the critical exponent for the semilinear Euler-Poisson-Darboux-Tricomi equation with power nonlinearity}

\author{Alessandro Palmieri}
\affil{Department of Mathematics, University of Bari, Via E. Orabona 4, 70125 Bari, Italy}

\maketitle

\begin{abstract}
In this note, we derive a blow-up result for a semilinear generalized Tricomi equation with damping and mass terms having time-dependent coefficients. We consider these coefficients with critical decay rates. Due to this threshold nature of the time-dependent coefficients (both for the damping and for the mass), the multiplicative constants appearing in these lower order terms strongly influence the value of the critical exponent, determining a competition between a Fujita-type exponent and a Strauss-type exponent. 
\end{abstract}

\begin{flushleft}
\textbf{Keywords} Critical exponent, Fujita exponent, Strauss exponent, Blow-up, Power nonlinearity
\end{flushleft}

\begin{flushleft}
\textbf{AMS Classification (2020)} 35B33, 35B44, 35L15, 35L71.
\end{flushleft}

\section{Introduction}

In the present note, we prove a blow-up result for local in time weak solutions to the following semilinear Cauchy problem with \emph{power nonlinearity} $|u|^p$
\begin{align}\label{semilinear EPDT}
\begin{cases} \partial_t^2u-t^{2\ell}\Delta u + \mu t^{-1} \partial_t u+\nu^2 t^{-2} u=| u|^p, &  x\in \mathbb{R}^n, \ t>1,\\
u(1,x)=\varepsilon u_0(x), & x\in \mathbb{R}^n, \\ \partial_t u(1,x)=\varepsilon u_1(x), & x\in \mathbb{R}^n,
\end{cases}
\end{align} where $\ell >-1$, $\mu,\nu^2$ are nonnegative real constants, $p>1$ and $\varepsilon$ is a positive constant describing the size of Cauchy data. We consider the case with initial data taken at the time $t=1$, nevertheless, the results that we are going to prove are valid for data taken at any initial time $t=t_0>0$. 

 For $\ell=0$ and $\nu^2=0$, the linearized equation associated with the equation in \eqref{semilinear EPDT} is known as \emph{Euler-Poisson-Darboux equation} (see the introduction of \cite{DAbb21} for a detailed overview on the literature regarding this model), while for $\mu=\nu^2=0$ the second-order operator $\partial^2_t-t^{2\ell}\Delta$ on the left-hand side of the equation in \eqref{semilinear EPDT} is called \emph{generalized Tricomi operator}. Motivated by these special cases and for the sake of brevity, we will call the equation in \eqref{semilinear EPDT} \emph{semilinear Euler-Poisson-Darboux-Tricomi equation} (semilinear EPDT equation). 

Over the last two decades, several papers have been devoted to the study of semilinear models with power nonlinearity, that are special cases of \eqref{semilinear EPDT} for given values of the parameters $\ell$, $\mu$, $\nu^2$.

Concerning the Cauchy problem associated with the semilinear generalized Tricomi equation
\begin{align}\label{semilinear Tricomi}
\begin{cases} \partial_t^2u-t^{2\ell}\Delta u =f(u,\partial_t u), &  x\in \mathbb{R}^n, \ t>0,\\
u(0,x)=\varepsilon u_0(x), & x\in \mathbb{R}^n, \\ \partial_t u(0,x)=\varepsilon u_1(x), & x\in \mathbb{R}^n,
\end{cases}
\end{align} we recall \cite{DL03,Yag04,Yag06} for the first results in the case with power nonlinearity $f(u,\partial_t u)=|u|^p$. More specifically, in \cite{DL03} some blow-up results are proved for weak solutions to the Cauchy problem associated with the generalized Tricomi operator (and, more generally, for Grushin-type operators) by mean of the test function method. The so-called \emph{quasi-homogeneous dimension} of the generalized Tricomi operator $\mathcal{Q}=(\ell+1)n+1$ made its appearance in the upper bound for $p$ in these results (see \cite{DL13,L18} for the definition of the quasi-homogeneous dimension for a more general partial differential operator).  On the other hand, in \cite{Yag04,Yag06} the fundamental solution for the generalized Tricomi operator (sometimes called also Gellerstedt operator in the literature) is employed to derive an integral representation formula, that is used in turn to prove (under suitable assumptions on $p$) the local existence of solutions and the global existence of small data solutions, respectively.

 Afterwards, in the series of papers \cite{HWY17,HWY16,HWY17d2,TL19,HWY20,Sun21} it was established that the critical exponent for \eqref{semilinear Tricomi} with $f(u,\partial_t u)=|u|^p$ and $n\geqslant 2$ is the Strauss-type exponent (cf. \cite{Str81} and the literature citing this renowned paper), that we denote $p_{\mathrm{Str}}(n,\ell)$ in the present paper, given by the biggest root of the quadratic equation
\begin{equation}\label{crit exponent tricomi}
\left(\frac{n-1}{2}+\frac{\ell}{2(\ell+1)}\right)p^2-\left(\frac{n+1}{2}-\frac{3\ell}{2(\ell+1)}\right)p-1=0.
\end{equation}
We recall that for us the fact that $p_{\mathrm{Str}}(n,\ell)$ is the \emph{critical exponent} for \eqref{semilinear Tricomi} with $f(u,\partial_t u)=|u|^p$ means the following: for any $1<p<p_{\mathrm{Str}}(n,\ell)$ local in time solutions blow up in finite time under suitable sign assumptions for the Cauchy data and regardless of their size, while for $p>p_{\mathrm{Str}}(n,\ell)$ (technical upper bounds for $p$ may appear, depending on the space for the solutions) a global in time existence result for small data solutions holds.

Very recently, even the cases with derivative type nonlinearity $f(u,\partial_t u)=|\partial_t u|^p$ and with mixed nonlinearity  $f(u,\partial_t u)=|u|^q+|\partial_t u|^p$ have been studied from the point of view of blow-up dynamics in \cite{LP20,CLP20,LS20,HH20Tri}. In particular, in \cite{LP20} a blow-up result when $f(u,\partial_t u)=|\partial_t u|^p$ is proved for $1<p\leqslant  \frac{\mathcal{Q}}{\mathcal{Q}-2}$.

Regarding the semilinear Euler-Poisson-Darboux equation and, more in general, the semilinear wave equation with scale-invariant damping and mass (i.e. \eqref{semilinear EPDT} for $\ell=0$), combining the contributions from many different authors (see \cite{DAbb15,DLR15,NPR16,IS17,PalRei18,PT18,DabbPal18,LST20,DAbb21,CGL21}
 and references therein) we may reasonably conjecture that for nonnegative $\mu, \nu^2$ such that the quantity
\begin{align}\label{def delta}
\delta\doteq (\mu-1)^2-4\nu^2
\end{align} satisfies $\delta\geqslant 0$ the critical exponent is given by
\begin{align} \label{critical exponent scale invariant}
\max\left\{p_{\mathrm{Str}}\left(n+\mu,0\right), \, p_{\mathrm{Fuj}}\left(n+\tfrac{\mu-1}{2}-\tfrac{\sqrt{\delta}}{2}\right) \right\}.
\end{align} Here $p_{\mathrm{Fuj}}(n)\doteq 1+\frac{2}{n}$ denotes the so-called \emph{Fujita exponent} (cf. \cite{Fuj66}) and its presence can be justified as a consequence of diffusion phenomena between the solutions to the corresponding linearized model and those to some suitable parabolic equation. We point out that the condition $\delta\geqslant 0$ implies somehow that the damping term $\mu t^{-1} \partial_t u$ has a dominant influence over the mass term $\nu^2 t^{-2} u$. Although the proof of the necessity part of this conjecture is fully demonstrated (see \cite{PalRei18,PT18,LST20}), for the sufficiency part only the one-dimensional case was recently completely clarified \cite{DAbb21}, while in the higher dimensional case only a few special cases have been studied (often in the radially symmetric case).

We emphasize that for $\ell=-\frac{2}{3}$ and $\mu=2,\nu^2=0$ the model in \eqref{semilinear EPDT} is the semilinear wave equation in the \emph{Einstein-de Sitter spacetime} with power nonlinearity (see \cite{GalYag17EdS}). Hence, for the semilinear wave equation in the generalized Einstein-de Sitter spacetime, i.e., when $\ell\in (-1,0)$ and $\mu\geqslant 0, \nu^2=0$ in \eqref{semilinear EPDT},  in the series of papers \cite{Pal20EdS,Pal20EdSmu,TW20,TW20h} several blow-up results are proved provided that $p>1$ is below  
\begin{align} \label{critical exponent EdS}
 \max\left\{ p_{\mathrm{Str}}\left(n+\tfrac{\mu}{\ell+1},\ell\right), \, p_{\mathrm{Fuj}}\left((\ell+1)n\right)\right\}.
 \end{align} 
Moreover, under the same conditions for the parameters as above, also the case with derivative type nonlinearity $|\partial_t u|^p$ is studied in \cite{HHP20,HHP20mu,TW21Gla}. Furthermore, we point out that in \cite{TW21,TW21Gla} even the case $\ell\leqslant -1$ with $\nu^2=0$  is studied for different semilinear terms.

The purpose of the present paper is twofold: on the one hand, we want to prove the necessity part for the one-dimensional case in \eqref{semilinear EPDT} that together with the sufficiency part from \cite{DAbb21} (cf. Corollary 5.1) will show the optimality of our result for $n=1$; on the other hand, we want to generalize the condition that   determines the critical exponent for \eqref{semilinear EPDT} in the general $n$-dimensional case, obtaining as a candidate to be critical an exponent that is consistent with all previously mentioned special cases.
Finally, we mention that in \cite{BHHT21} the semilinear Cauchy problem with derivative type nonlinearity and the same linear partial differential operator as in \eqref{semilinear EPDT} is considered, and a Glassey-type exponent is found as the upper bound for the exponent in the blow-up range.

\paragraph{Notations} Throughout the paper we denote by
\begin{align}\label{def phi ell}
\phi_\ell(t)\doteq \frac{t^{\ell+1}}{\ell+1} 
\end{align} the primitive of the speed of propagation $t^{\ell}$ that vanishes for $t=0$. In particular, the amplitude of the light-cone for the Cauchy problem with data prescribed at the initial time $t_0=1$ is given by the function $\phi_\ell(t)-\phi_\ell(1)$. 
The ball in $\mathbb{R}^n$ with radius $R$ around the origin is denoted $B_R$. The notation $f\lesssim g$ means that there exists a positive constant $C$ such that $f\leqslant Cg$ and, analogously, $f\gtrsim g$.
Finally, as in the introduction, we will denote by $p_{\mathrm{Fuj}}(n)$ the Fujita exponent and by $p_{\mathrm{Str}}(n,\ell)$ the Strauss-type exponent. 

\subsection{Main results}

Let us begin this section by introducing the notion of weak solutions to \eqref{semilinear EPDT} that we will employ throughout the entire paper.

\begin{definition} \label{def sol} Let $u_0,u_1\in L^1_{\mathrm{loc}}(\mathbb{R}^n)$ such that $\mathrm{supp} \, u_0, \mathrm{supp} \, u_1 \subset B_R$ for some $R>0$. We say that $$u\in\mathcal{C}\big([1,T),W^{1,1}_{\mathrm{loc}}(\mathbb{R}^n)\big)\cap \mathcal{C}^1\big([1,T),L^{1}_{\mathrm{loc}}(\mathbb{R}^n)\big)\cap L^p_{\mathrm{loc}}\big((1,T)\times \mathbb{R}^n\big)$$ is a \emph{weak solution} to \eqref{semilinear EPDT} on $[1,T)$ if $u(1,\cdot)=\varepsilon u_0$ in $L^{1}_{\mathrm{loc}}(\mathbb{R}^n)$, $u$ fulfills the support condition 
\begin{align}\label{support condition sol}
\mathrm{supp} \, u(t,\cdot) \subset B_{R+\phi_\ell(t)-\phi_\ell(1)} \qquad \mbox{for any} \ t\in(1,T),
\end{align} 
 and the integral identity
\begin{align}
& \int_{\mathbb{R}^n} \partial_t u(t,x)\phi(t,x) \, \mathrm{d}x+\int_1^t\int_{\mathbb{R}^n} \big(-\partial_t u(s,x)\phi_s(s,x)+s^{2\ell}\nabla u(s,x) \cdot \nabla \phi(s,x)\big)\mathrm{d}x\, \mathrm{d}s \notag  \\
 & \qquad + \int_1^t\int_{\mathbb{R}^n} \big(\mu s^{-1}\partial_t u(s,x)\phi(s,x)+\nu^2 s^{-2} u(s,x)  \phi(s,x)\big)\mathrm{d}x\, \mathrm{d}s \notag \\
& \quad = \varepsilon \int_{\mathbb{R}^n} u_1(x)\phi(1,x) \, \mathrm{d}x+ \int_1^t\int_{\mathbb{R}^n}|u(s,x)|^p\phi(s,x) \, \mathrm{d}x\, \mathrm{d}s \label{integral rel def sol}
\end{align} holds for any $t\in(1,T)$ and any test function $\phi\in\mathcal{C}^\infty_0\big([1,T)\times\mathbb{R}^n\big)$.
\end{definition}

We notice that, performing further steps of integration by parts in \eqref{integral rel def sol}, we get the integral relation
\begin{align}
& \int_{\mathbb{R}^n} \big( \partial_t u(t,x)\phi(t,x)- u(t,x)\phi_t(t,x)+\mu t^{-1} u(t,x)\phi(t,x) \big) \mathrm{d}x \notag  \\
 & \qquad +\int_1^t\int_{\mathbb{R}^n} u(s,x) \big( \phi_{ss}(s,x)-s^{2\ell}\Delta \phi(s,x)-\mu s^{-1}\phi_s(s,x)+(\mu+\nu^2) s^{-2}  \phi(s,x)\big)\mathrm{d}x\, \mathrm{d}s \notag \\
& \quad = \varepsilon \int_{\mathbb{R}^n} \big( (u_1(x)+\mu u_0(x))\phi(1,x) -u_0(x) \phi_t(1,x)\big) \mathrm{d}x+ \int_1^t\int_{\mathbb{R}^n}|u(s,x)|^p\phi(s,x) \, \mathrm{d}x\, \mathrm{d}s \label{integral rel def sol 2}
\end{align} for any $t\in(1,T)$ and any test function $\phi\in\mathcal{C}^\infty_0\big([1,T)\times\mathbb{R}^n\big)$.

Let us state our result in the sub-critical case.
\begin{theorem}\label{Thm main} Let $\ell>-1$ and $\mu,\nu^2\geqslant 0$ such that $\delta\geqslant 0$. Let us assume that the exponent $p$ of the nonlinear term satisfies
\begin{align*}
1<p< \max\left\{p_{\mathrm{Str}}\!\left(n+\tfrac{\mu}{\ell+1},\ell\right),p_{\mathrm{Fuj}}\!\left((\ell+1)n+\tfrac{\mu-1}{2}-\tfrac{\sqrt{\delta}}{2}\right) \right\}.
\end{align*} 
Let $u_0,u_1\in L^1_{\mathrm{loc}}(\mathbb{R}^n)$ be nonnegative,  nontrivial and compactly supported functions with supports contained in  $B_R$ for some $R>0$ such that
\begin{align}\label{integral assumption Cauchy data}
u_1 +\tfrac{\mu-1-\sqrt{\delta}}{2} u_0 \geqslant 0. 
\end{align} 
Let $u\in\mathcal{C}\big([1,T),W^{1,1}_{\mathrm{loc}}(\mathbb{R}^n)\big)\cap \mathcal{C}^1\big([1,T),L^{1}_{\mathrm{loc}}(\mathbb{R}^n)\big)\cap L^p_{\mathrm{loc}}\big((1,T)\times \mathbb{R}^n\big)$ be a weak solution to \eqref{semilinear EPDT} according to Definition \ref{def sol} with lifespan $T=T(\varepsilon)$.

Then, there exists a positive constant $\varepsilon_0 = \varepsilon_0(u_0,u_1,n,p,\ell,\mu,\nu^2,R)$ such that for any $\varepsilon \in (0,\varepsilon_0]$ the weak solution $u$ blows up in finite time. Furthermore, the upper bound estimates for the lifespan
\begin{align}\label{lifespan est thm main}
T(\varepsilon)\leqslant \begin{cases} 
C\varepsilon^{-\frac{p(p-1)}{\theta(n,\ell,\mu,p)}} & \mbox{if} \ p<p_{\mathrm{Str}}\big(n+\tfrac{\mu}{\ell+1},\ell\big), \\
C\varepsilon^{-\big(\frac{2}{p-1}-\big((\ell+1)n+\frac{\mu-1}{2}-\frac{\sqrt{\delta}}{2}\big)\big)^{-1}} & \mbox{if} \ p<p_{\mathrm{Fuj}}\big((\ell+1)n+\tfrac{\mu-1}{2}-\tfrac{\sqrt{\delta}}{2}\big),
\end{cases}
\end{align}
holds, where the positive constant $C$ is independent of $\varepsilon$ and 
\begin{align}\label{def theta}
\theta(n,\ell,\mu,p)\doteq \ell+1 +\Big(\tfrac{n+1}{2}(\ell+1)+\tfrac{\mu-3\ell}{2}\Big) p-\Big(\tfrac{n-1}{2}(\ell+1)+\tfrac{\ell+\mu}{2}\Big) p^2.
\end{align}
\end{theorem}

\begin{remark} In the previous statement it might happen that the argument of the Fujita exponent is a nonpositive number (however, only for $\ell<0$ and $\mu \in [0,1)$). Whenever this happens, we do not require any upper bound for $p>1$, meaning formally that $p_{\mathrm{Fuj}}(k)=\infty$ for $k\leqslant 0$.
\end{remark}

\begin{remark} The exponent $p_{\mathrm{Str}}\big(n+\tfrac{\mu}{\ell+1},\ell\big)$ is obtained from the Strauss-type exponent $p_{\mathrm{Str}}(n,\ell)$ defined through \eqref{crit exponent tricomi} by a shift of magnitude $\tfrac{\mu}{\ell+1}$ in the space dimension. Equivalently, $p_{\mathrm{Str}}\big(n+\tfrac{\mu}{\ell+1},\ell\big)$ is the positive root to the quadratic equation
\begin{align*}
\bigg(\frac{n-1}{2}(\ell+1)+\frac{\ell+\mu}{2}\bigg) p^2 -\bigg(\frac{n+1}{2}(\ell+1)+\frac{\mu-3\ell}{2}\bigg) p -(\ell+1) =0.
\end{align*}
\end{remark}

Finally, we provide a blow-up result when we consider the critical Fujita-type exponent.

\begin{theorem}\label{Thm crit Fuj} Let $\ell>-1$ and $\mu,\nu^2\geqslant 0$ such that $\delta\geqslant 0$. Let us assume that the exponent $p$ of the nonlinear term satisfies
\begin{align*}
p=p_{\mathrm{Fuj}}\left((\ell+1)n+\tfrac{\mu-1}{2}-\tfrac{\sqrt{\delta}}{2}\right) .
\end{align*} 
Let $u_0,u_1\in L^1_{\mathrm{loc}}(\mathbb{R}^n)$ be nonnegative,  nontrivial and compactly supported functions with supports contained in  $B_R$ for some $R>0$.
Let $u\in\mathcal{C}\big([1,T),W^{1,1}_{\mathrm{loc}}(\mathbb{R}^n)\big)\cap \mathcal{C}^1\big([1,T),L^{1}_{\mathrm{loc}}(\mathbb{R}^n)\big)\cap L^p_{\mathrm{loc}}\big((1,T)\times \mathbb{R}^n\big)$ be a weak solution to \eqref{semilinear EPDT} according to Definition \ref{def sol} with lifespan $T=T(\varepsilon)$.

Then, there exists a positive constant $\varepsilon_0 = \varepsilon_0(u_0,u_1,n,p,\ell,\mu,\nu^2,R)$ such that for any $\varepsilon \in (0,\varepsilon_0]$ the weak solution $u$ blows up in finite time. Moreover, the upper bound estimates for the lifespan
\begin{align}\label{lifespan est thm crit Fuj}
T(\varepsilon)\leqslant \begin{cases} 
\exp\left(E\varepsilon^{-(p-1)}\right) & \mbox{if} \ \delta>0, \\
\exp\left(E\varepsilon^{-(p-1)/p}\right) & \mbox{if} \ \delta=0,
\end{cases}
\end{align}
holds, where the positive constant $E$ is independent of $\varepsilon$.
\end{theorem}

\section{Proof of Theorem \ref{Thm main}}

In this section, we prove Theorem \ref{Thm main} by deriving a sequence of lower bound estimates for the space average of a weak solution $u$ to \eqref{semilinear EPDT}. More precisely, introducing the functional
\begin{align*}
U(t)\doteq \int_{\mathbb{R}^n} u(t,x)\, \mathrm{d}x \qquad \mbox{for} \ t\in [1,T) ,
\end{align*} our aim is to determine estimates from below for $U$.
 Letting to $\infty$ the index of this sequence of lower bounds, we establish that the space average of $u$ cannot be globally in time defined and we determine as a byproduct an upper bound estimate for the lifespan. In particular, the first two steps that we need to carry out in order to apply the iteration argument are determining the \emph{iteration frame} for $U$ (namely, an integral inequality, where $U$ appears both on the left and right-hand side) and a first lower bound estimate for $U$. In order to establish such a first lower bound estimate for $U$ we consider a suitable positive solution to the adjoint equation to the linear EPDT equation. Finally, we employ this first lower bound estimate for $U$ to begin the iteration procedure, plugging it in the iteration frame. Then, repeating iteratively the procedure we determine the desired sequence of lower bounds.

\subsection{Derivation of the iteration frame}

In this subsection we derive the iteration frame for $U$. To this purpose we employ the \emph{double multiplier technique} from \cite{PT18} (see also \cite{ITW21,PalTak22}). Given $t\in (1,T)$, let us begin by choosing as test function in \eqref{integral rel def sol} a cut-off function that localizes the forward light-cone, namely, we take $\phi\in\mathcal{C}^\infty_0\big([1,T)\times \mathbb{R}^n\big)$ such that $\phi=1$ in $\{(s,x)\in [1,t]\times \mathbb{R}^n: |x|\leqslant R+\phi_\ell(s)-\phi_\ell(1)\}$. Therefore, from \eqref{integral rel def sol} we get
\begin{align*}
& \int_{\mathbb{R}^n} \partial_t u(t,x) \, \mathrm{d}x +\int_1^t \int_{\mathbb{R}^n} \big(\mu s^{-1}\partial_t u(s,x)+\nu^2 s^{-2}u(s,x)\big) \mathrm{d}x \, \mathrm{d}s \\ & \qquad = \varepsilon \int_{\mathbb{R}^n} u_1(x) \, \mathrm{d}x +\int_1^t \int_{\mathbb{R}^n} |u(s,x)|^p \,  \mathrm{d}x \, \mathrm{d}s.
\end{align*} Differentiating with respect to $t$ the previous equality, we obtain
\begin{align}
\int_{\mathbb{R}^n} |u(t,x)|^p \,  \mathrm{d}x & = \int_{\mathbb{R}^n} \partial_t^2 u(t,x) \, \mathrm{d}x + \mu t^{-1} \int_{\mathbb{R}^n} \partial_t u(t,x) \, \mathrm{d}x  + \nu^2 t^{-2}  \int_{\mathbb{R}^n}u(t,x) \, \mathrm{d}x \notag  \\
& = U''(t) + \mu t^{-1} U'(t)+ \nu^2 t^{-2} U(t). \label{ODE for U}
\end{align} If we denote by $r_1,r_2$ the roots of the quadratic equation $$r^2-(\mu-1) r+\nu^2=0,$$ then, we can rewrite the right-hand side of the last relation as follows
\begin{align} \label{double multiplier}
 U''(t) + \mu t^{-1} U'(t)+ \nu^2 t^{-2} U(t) =t^{-(r_2+1)}\frac{\mathrm{d}}{\mathrm{d}t}\Big( t^{r_2+1-r_1}\frac{\mathrm{d}}{\mathrm{d}t}\Big(t^{r_1}U(t)\Big)\!\Big).
\end{align} We emphasize that the role of $r_1$ and $r_2$ is fully interchangeable in the previous identity.
Combining \eqref{ODE for U} and \eqref{double multiplier}, after some straightforward steps we arrive at
\begin{align} \label{OIE for U}
U(t) &=U_{\mathrm{lin}}(t) + \int_1^t \left(\frac{s}{t}\right)^{r_1} \int_1^s \left(\frac{\tau}{s}\right)^{r_2+1}\int_{\mathbb{R}^n}|u(\tau,x)|^p \, \mathrm{d}x \, \mathrm{d}\tau \, \mathrm{d}s,
\end{align} 
where 
\begin{align}\label{U lin def}
U_{\mathrm{lin}}(t) \doteq \begin{cases}
\displaystyle{\tfrac{r_1 t^{-r_2}-r_2 t^{-r_1}}{r_1-r_2}\, U(1)+\tfrac{t^{-r_2}-t^{-r_1}}{r_1-r_2} \, U'(1)} & \mbox{if} \ \ \delta>0, \\
t^{-r_1}(1+r_1\ln t)\, U(1)+t^{-r_1}\ln t \, U'(1) & \mbox{if} \ \ \delta=0.
\end{cases}
\end{align}
Let us set $$r_1\doteq \tfrac{\mu-1-\sqrt{\delta}}{2}, \quad r_2\doteq \tfrac{\mu-1+\sqrt{\delta}}{2},$$ where the definition of $\delta$ is given in \eqref{def delta}. From here on, these will be the fixed values of $r_1$, $r_2$.
Consequently, from \eqref{OIE for U} we have a twofold result. On the one hand, we get the lower bound estimate for $U$
\begin{align}\label{lower bound U Fujita type}
U(t)\geqslant  I \varepsilon \, t^{-r_1} \qquad \mbox{for} \ t \in [1,T),
\end{align} where the multiplicative constant $I$ depends on the positive quantities $\int_{\mathbb{R}^n} u_0(x)\, \mathrm{d}x$ and $\int_{\mathbb{R}^n} u_1(x)\, \mathrm{d}x$. 

On the other hand, by using again 
the nonnegativity and nontriviality of $u_0,u_1$, we find 
\begin{align}
U(t) &=  \int_1^t \left(\frac{s}{t}\right)^{r_1} \int_1^s \left(\frac{\tau}{s}\right)^{r_2+1}\int_{\mathbb{R}^n}|u(\tau,x)|^p \, \mathrm{d}x \, \mathrm{d}\tau \, \mathrm{d}s   \label{pre iteration frame for U} \\
& \gtrsim  \int_1^t \left(\frac{s}{t}\right)^{r_1} \int_1^s \left(\frac{\tau}{s}\right)^{r_2+1} (R+\phi_\ell(\tau)-\phi_\ell(1))^{-n(p-1)} (U(\tau))^p\, \mathrm{d}\tau \, \mathrm{d}s  \notag\\
& \gtrsim   t^{-r_1} \int_1^t s^{r_1-r_2-1} \int_1^s \tau^{r_2+1} (1+\tau)^{-n(\ell+1)(p-1)} (U(\tau))^p\, \mathrm{d}\tau \, \mathrm{d}s. \notag 
\end{align}  Note that in the last chain of inequalities we used the support condition for $u$ and Jensen's inequality.
Hence, we obtained the following the iteration frame for $U$
\begin{align}
U(t) & \geqslant C t^{-r_1} \int_1^t s^{r_1-r_2-1} \int_1^s \tau^{r_2+1} \tau^{-n(\ell+1)(p-1)} (U(\tau))^p\, \mathrm{d}\tau \, \mathrm{d}s 
\label{iteration frame for U}
\end{align} for a suitable positive constant $C$ that depends on $n,p,\ell$.

\subsection{Solution of the adjoint equation}
In the previous subsection we established a first lower bound estimate for $U$ in \eqref{lower bound U Fujita type}. This estimate will be the starting point for the proof of the blow-up result for $p<p_{\mathrm{Fuj}}\big((\ell+1)n+\tfrac{\mu-1}{2}-\tfrac{\sqrt{\delta}}{2}\big)$. However, to prove the blow-up result for $p<p_{\mathrm{Str}}\big(n+\tfrac{\mu}{\ell+1},\ell\big)$ we need to determine a further lower bound estimate for $U$. According to this purpose, we introduce a second auxiliary time-dependent functional, which is a certain weighted spatial average of $u$. In particular, we are going to choose the weight function to be a positive solution of the adjoint equation to the linear EPDT equation, namely,
\begin{align}\label{adj eq EPDT}
\partial_s^2 \psi -s^{2\ell}\Delta \psi-\tfrac{\partial}{\partial s}(\mu s^{-1} \psi)+\nu^2 s^{-2} \psi=0.
\end{align} 
In order to determine a suitable solution to \eqref{adj eq EPDT}, we follow the approach from \cite[Subsection 2.1]{Pal20EdSmu}. We look for a solution to \eqref{adj eq EPDT} with separate variables. As $x$-dependent function we consider the function 
\begin{align*}
\varphi(x)\doteq \begin{cases} \mathrm{e}^{x}+\mathrm{e}^{- x} & \mbox{if} \ \, n=1, \\ \int_{\mathbb{S}^{n-1}}\mathrm{e}^{x\cdot \omega} \mathrm{d} \sigma_\omega & \mbox{if} \ \, n\geqslant 2,   \end{cases}
\end{align*} that has been introduce for the study of blow-up phenomena for semilinear hyperbolic models in \cite{YZ06}. The positive function $\varphi\in\mathcal{C}^\infty(\mathbb{R}^n)$ satisfies $\Delta \varphi=\varphi$ and
\begin{align}\label{asymptotic phi}
\varphi(x) \sim c_n |x|^{-\frac{n-1}{2}} \mathrm{e}^{|x|} \qquad \mbox{as} \ |x|\to \infty,
\end{align}  where $c_n$ is a positive constant depending on $n$.

On the other hand, as $s$-dependent function we look for a positive solution to the second-order linear ODE
\begin{align} \label{eq varrho wrt s}
\varrho''(s)-s^{2\ell}\varrho(s)-\mu s^{-1} \varrho'(s)+(\mu+\nu^2)s^{-2}\varrho(s)=0.
\end{align} Let us perform the change of variables $\sigma = \phi_\ell(s)$. Then, $\varrho$ solves the previous equation if and only if 
\begin{align}\label{eq varrho wrt sigma}
\sigma^2\frac{\mathrm{d}^2\varrho}{\mathrm{d}\sigma^2} +\frac{\ell-\mu}{\ell+1} \, \sigma \frac{\mathrm{d}\varrho}{\mathrm{d}\sigma} +\bigg(\frac{\mu+\nu^2}{(\ell+1)^2}-\sigma^2\bigg) \varrho=0.
\end{align} Next, we consider the transformation $\varrho(\sigma) = \sigma^\alpha \eta(\sigma)$ with $\alpha\doteq \frac{\mu+1}{2(\ell+1)}$. Hence, $\varrho$ solves \eqref{eq varrho wrt sigma} if and only if $\eta$ is a solution to 
\begin{align} \label{eq eta}
\sigma^2\frac{\mathrm{d}^2\eta}{\mathrm{d}\sigma^2} + \sigma \frac{\mathrm{d}\eta}{\mathrm{d}\sigma} -\Big(\sigma^2+ \frac{\delta}{4(\ell+1)^2}\Big) \eta=0.
\end{align}  As solution to \eqref{eq eta} we choose the modified Bessel function of the second kind $\mathrm{K}_{\frac{\sqrt{\delta}}{2(\ell+1)}}(\sigma)$. Consequently, we set as positive solutions to \eqref{eq varrho wrt s} 
\begin{align}\label{def varrho}
\varrho(s) &\doteq s^{\frac{\mu+1}{2}}\mathrm{K}_{\frac{\sqrt{\delta}}{2(\ell+1)}}(\phi_\ell(s)).
\end{align} Note that $\varrho=\varrho(s;\ell,\mu,\nu^2)$, but for the sake of brevity we will skip the dependence on $\ell,\mu,\nu^2$ in the notations hereafter.
Therefore, we may define now the following function as positive solution to \eqref{adj eq EPDT}
\begin{align}
\psi(s,x) &\doteq \varrho(s) \varphi(x). \label{def psi}
\end{align}
By using the weight function $\psi$, we introduce the auxiliary functional 
\begin{align*}
U_0(t)\doteq \int_{\mathbb{R}^n} u(t,x)\psi(t,x) \, \mathrm{d}x.
\end{align*} Notice that thanks to the support condition for the solution $u$ from Definition \ref{def sol}, it is possible to employ $\psi$ as test function in \eqref{integral rel def sol 2}. Consequently, using \eqref{adj eq EPDT}, we get
\begin{align}
& \int_{\mathbb{R}^n} \big( \partial_t u(t,x)\psi(t,x)- u(t,x)\psi_t(t,x)+\mu t^{-1} u(t,x)\psi(t,x) \big) \mathrm{d}x \notag  \\
& \quad = \varepsilon \int_{\mathbb{R}^n} \big( \varrho(1) (u_1(x)+\mu u_0(x))-  \varrho'(1) u_0(x)\big)\varphi(x) \mathrm{d}x+ \int_1^t\int_{\mathbb{R}^n}|u(s,x)|^p\psi(s,x) \, \mathrm{d}x\, \mathrm{d}s. \label{fundamental integral identity}
\end{align} Applying the recursive relation for the derivative of the modified Bessel function of the second kind $\partial_z \mathrm{K}_\gamma(z)=-\mathrm{K}_{\gamma+1}(z)+\frac{\gamma}{z}\mathrm{K}_\gamma(z)$ (cf. \cite[Equations (10.29.1)]{OLBC10}), we have
\begin{align*}
\varrho'(s)& =-s^{\frac{\mu+1}{2}+\ell} \mathrm{K}_{\frac{\sqrt{\delta}}{2(\ell+1)}+1}(\phi_\ell(s)) +\tfrac{\mu+1+\sqrt{\delta}}{2} \, s^{\frac{\mu-1}{2}} \mathrm{K}_{\frac{\sqrt{\delta}}{2(\ell+1)}}(\phi_\ell(s)).
\end{align*}  Thus,
\begin{align*}
I_{\ell,\mu,\nu^2}[u_0,u_1] & \doteq   \int_{\mathbb{R}^n} \big( \varrho(1) u_1(x)+(\varrho(1) \mu -  \varrho'(1)) u_0(x) \big)\varphi(x) \mathrm{d}x \\ & = \int_{\mathbb{R}^n} \Big(\mathrm{K}_{\frac{\sqrt{\delta}}{2(\ell+1)}+1}(\phi_\ell(1)) \,  u_0(x) +  \mathrm{K}_{\frac{\sqrt{\delta}}{2(\ell+1)}}(\phi_\ell(1)) \big(u_1(x)+\tfrac{\mu-1-\sqrt{\delta}}{2} u_0(x) \big)\Big)\varphi(x) \mathrm{d}x >0,
\end{align*} where we employed the nonnegativity and the nontriviality of $u_0$ and \eqref{integral assumption Cauchy data}.

Hence, we can rewrite \eqref{fundamental integral identity} as 
\begin{align*}
\varepsilon \, I_{\ell,\mu,\nu^2}[u_0,u_1] + \int_1^t\int_{\mathbb{R}^n}|u(s,x)|^p\psi(s,x) \, \mathrm{d}x\, \mathrm{d}s & = U'_0(t)+\mu t^{-1}U_0(t)  -\frac{2\varrho'(t)}{\varrho(t)} U_0(t) \\
&= \frac{\varrho^2(t)}{t^{\mu}}\frac{\mathrm{d}}{\mathrm{d}t} \left(\frac{t^\mu}{\varrho^2(t)} U_0(t)\right).
\end{align*} Since both $\psi$ and the nonlinear term are nonnegative, from the previous relation we obtain
\begin{align*}
U_0(t)\geqslant  \varepsilon \, I_{\ell,\mu,\nu^2}[u_0,u_1] \, \frac{\varrho^2(t)}{t^{\mu}}\int_1^t  \frac{s^{\mu}}{\varrho^2(s)}\,  \mathrm{d}s >0 \qquad \mbox{for any}  \ t\in [1,T).
\end{align*} Thanks to the assumptions on the Cauchy data we have shown that the functional $U_0$ is nonnegative. Next, we shall determine a lower bound estimate for $U_0$ for large times. For this reason we recall the asymptotic behavior of $\mathrm{K}_\gamma$ for large arguments, namely, $\mathrm{K}_\gamma(z)= \sqrt{\pi/(2z)}\mathrm{e}^{-z} (1+\mathcal{O}(z^{-1}))$ as $z\to \infty$ for $z>0$ (cf. \cite[Equation (10.25.3)]{OLBC10}). So, there exists $T_0=T_0(\ell,\mu,\nu^2)>1$ such that for $s\geqslant T_0$ it holds
\begin{align}\label{asymp varrho}
\tfrac{1}{4} \pi (\ell+1) \, \mathrm{e}^{-2\phi_\ell(s)} s^{\mu-\ell}\leqslant \varrho^2(s)\leqslant \pi (\ell+1)\,  \mathrm{e}^{-2\phi_\ell(s)} s^{\mu-\ell}.
\end{align} Then, for $t\geqslant T_0$ we have
\begin{align*}
U_0(t) & \geqslant  \varepsilon \, I_{\ell,\mu,\nu^2}[u_0,u_1] \, \frac{\varrho^2(t)}{t^{\mu}}\int_{T_0}^t  \frac{s^{\mu}}{\varrho^2(s)}\,  \mathrm{d}s \\
& \geqslant  \tfrac{\varepsilon}{4} \, I_{\ell,\mu,\nu^2}[u_0,u_1] \, t^{-\ell}  \mathrm{e}^{-2\phi_\ell(t)} \int_{T_0}^t  s^{\ell}  \mathrm{e}^{2\phi_\ell(s)} \,  \mathrm{d}s.
\end{align*} For $t\geqslant  2T_0$, it results
\begin{align*}
U_0(t) &  \gtrsim \varepsilon \, I_{\ell,\mu,\nu^2}[u_0,u_1] \, t^{-\ell}  \mathrm{e}^{-2\phi_\ell(t)} \int_{t/2}^t  s^{\ell}  \mathrm{e}^{2\phi_\ell(s)} \,  \mathrm{d}s \\
&  \gtrsim \varepsilon \, I_{\ell,\mu,\nu^2}[u_0,u_1] \, t^{-\ell}  \big(1- \mathrm{e}^{2\phi_\ell(t/2)-2\phi_\ell(t)}\big) = \varepsilon \, I_{\ell,\mu,\nu^2}[u_0,u_1] \, t^{-\ell}  \Big(1- \mathrm{e}^{-\frac{2}{\ell+1}(2^{\ell+1}-1)t^{\ell+1}}\Big) \\
& \gtrsim \varepsilon \, I_{\ell,\mu,\nu^2}[u_0,u_1] \, t^{-\ell} .
\end{align*} We emphasize that in the last step we use the condition $\ell>-1$ to estimate from below the factor containing the exponential term with a positive constant.

 Summarizing, we proved
\begin{align}\label{lower bound U0}
U_0(t)\gtrsim \varepsilon t^{-\ell} \qquad \mbox{for} \ t\geqslant 2T_0,
\end{align} where the unexpressed multiplicative constant depends on the Cauchy data and on the parameters $\ell,\mu,\nu^2$. Finally, we show how from \eqref{lower bound U0} we derive a second lower bound estimate for $U$.
By H\"older's inequality, we find  for $t\geqslant 2T_0$
\begin{align*}
 \varepsilon t^{-\ell} \lesssim U_0(t) \leqslant \left(\,\int_{\mathbb{R}^n}|u(t,x)|^p \, \mathrm{d}x\right)^{1/p} \left(\int_{B_{R+\phi_\ell(t)-\phi_\ell(1)}} (\psi(t,x))^{p'}\, \mathrm{d}x \right)^{1/p'},
\end{align*} where $p'$ denotes the conjugate exponent of $p$.
Following \cite[Section 3]{PalRei18} and employing \eqref{asymp varrho}, we can estimate for $t\geqslant 2T_0$
\begin{align*}
\int_{\mathbb{R}^n}|u(t,x)|^p \, \mathrm{d}x & \gtrsim \varepsilon^p t^{-\ell p} \left(\int_{B_{R+\phi_\ell(t)-\phi_\ell(1)}} (\psi(t,x))^{p'}\, \mathrm{d}x \right)^{-(p-1)} \\
& \gtrsim \varepsilon^p t^{-\ell p} (\varrho(t))^{-p} \mathrm{e}^{-p(R+\phi_\ell(t)-\phi_\ell(1))} (R+\phi_\ell(t)-\phi_\ell(1))^{-(n-1)(p-1)+\frac{n-1}{2}p} \\
& \gtrsim \varepsilon^p t^{\left(-\frac{n-1}{2}(\ell+1)-\frac{\ell+\mu}{2}\right)p+(n-1)(\ell+1)}.
\end{align*} Plugging this last estimate from below for the space integral of the nonlinear term in \eqref{pre iteration frame for U}, for any $t\geqslant T_1\doteq 2T_0+1$ we get
\begin{align}\label{1st lower bound U Str}
U(t)\geqslant K \varepsilon^p t^{\left(-\frac{n-1}{2}(\ell+1)-\frac{\ell+\mu}{2}\right)p+(n-1)(\ell+1)+2},
\end{align} where the multiplicative constant $K$ depends on the Cauchy data and on $n,p,\ell,\mu,\nu^2,R$, which is the desired lower bound estimate for $U$.

\subsection{Iteration argument}

In this section we establish the following sequence of lower bound estimates for $U$
\begin{align} \label{lowe bound U j}
U(t)\geqslant C_j t^{-\alpha_j} (t-T_1)^{\beta_j} \quad \mbox{for} \ t\geqslant T_1,
\end{align} where $\{\alpha_j\}_{j\in\mathbb{N}}$, $\{\beta_j\}_{j\in\mathbb{N}}$ and $\{C_j\}_{j\in\mathbb{N}}$ are sequences of positive real numbers that we will determine iteratively during the proof.

Let us begin by considering the case $1<p<p_{\mathrm{Str}}\big(n+\frac{\mu}{\ell+1},\ell\big)$. As we have previously mentioned, we consider \eqref{1st lower bound U Str} as first lower bound estimate for $U$ for $p$ below the Strauss-type exponent. Therefore, \eqref{lowe bound U j} for $j=0$ is satisfied provided that we set
$C_0\doteq K\varepsilon^p$, $\alpha_0\doteq [\frac{n-1}{2}(\ell+1)+\frac{\ell+\mu}{2}]p$ and $\beta_0\doteq (n-1)(\ell+1)+2$. 

Next, we assume that \eqref{lowe bound U j} holds true for some $j\geqslant 0$ and we want to prove it for $j+1$. If we plug the lower bound estimate for $U$ in \eqref{lowe bound U j} in the iteration frame \eqref{iteration frame for U}, we get 
\begin{align*}
U(t) & \geqslant C   t^{-r_1} \int_{T_1}^t s^{r_1-r_2-1} \int_{T_1}^s \tau^{r_2+1-n(\ell+1)(p-1)} (U(\tau))^p\, \mathrm{d}\tau \, \mathrm{d}s \\
& \geqslant C C_j^p    t^{-r_1} \int_{T_1}^t s^{r_1-r_2-1} \int_{T_1}^s (\tau-T_1)^{r_2+1+p\beta_j} \tau^{-n(\ell+1)(p-1)-p \alpha_j} \, \mathrm{d}\tau \, \mathrm{d}s \\
& \geqslant C C_j^p    t^{-r_2-1-n(\ell+1)(p-1)-p \alpha_j} \int_{T_1}^t \int_{T_1}^s (\tau-T_1)^{r_2+1+p\beta_j} \, \mathrm{d}\tau \, \mathrm{d}s \\
& \geqslant \frac{C C_j^p}{ (r_2+3+p\beta_j)^2}  \,  t^{-r_2-1-n(\ell+1)(p-1)-p \alpha_j}  (t-T_1)^{r_2+3+p\beta_j} ,
\end{align*} where we used the relations $r_1-r_2-1=-(\sqrt{\delta}+1)<0$ and $r_2+1>0$. Hence, if we define
\begin{align}
C_{j+1} & \doteq C C_j^p(r_2+3+p\beta_j)^{-2}, \label{def Cj}\\
 \alpha_{j+1} &\doteq r_2+1+n(\ell+1)(p-1)+p \alpha_j, \quad \beta_{j+1} \doteq r_2+3+p \beta_j, \label{def alpha j beta j}
\end{align}
then, we proved \eqref{lowe bound U j} for $j+1$ as well.
Let us determine explicitly $\alpha_j$ and $\beta_j$. By using recursively the definition in \eqref{def alpha j beta j}, we find
\begin{align}
\alpha_j &=  r_2+1+n(\ell+1)(p-1)+p \alpha_{j-1} = \ldots = ( r_2+1+n(\ell+1)(p-1)) \sum_{k=0}^{j-1} p^k + p^j \alpha_0 \notag \\ &= \left(\frac{r_2+1}{p-1}+n(\ell+1)+\alpha_0\right)p^j-\frac{r_2+1}{p-1}-n(\ell+1), \label{alpha j}
\end{align} and, analogously,
\begin{align}\label{beta j}
\beta_j = \left(\frac{r_2+3}{p-1}+\beta_0\right)p^j-\frac{r_2+3}{p-1}.
\end{align} In particular, combining \eqref{def Cj} and \eqref{beta j}, we get
\begin{align*}
C_{j} & = C C_{j-1}^p(r_2+3+p\beta_{j-1})^{-2}=  C C_{j-1}^p \beta_j^{-2} \geqslant  \underbrace{C \left(\frac{r_2+3}{p-1}+\beta_0\right)^{-2}}_{\doteq D} C_{j-1}^p p^{-2j} = D C_{j-1}^p p^{-2j} 
\end{align*} for any $j\geqslant 1$. Applying the logarithmic function to both sides of the inequality $C_{j}\geqslant  D C_{j-1}^p p^{-2j} $ and employing iteratively the resulting relation, we have
\begin{align*}
\ln C_j & \geqslant p \ln C_{j-1}-2j\ln p+\ln D \\
& \geqslant p^2 \ln C_{j-2}-2(j+(j-1)p)\ln p+(1+p) \ln D  \\ & \geqslant \ldots \geqslant p^j \ln C_0 - 2\ln p\left( \, \sum_{k=0}^{j-1} (j-k)p^k\right)+\ln D \sum_{k=0}^{j-1}p^k \\
&= p^j\left(\ln C_0-\frac{2p\ln p}{(p-1)^2}+\frac{\ln D}{p-1}\right)+\frac{2\ln p}{p-1}\, j+\frac{2p\ln p}{(p-1)^2}-\frac{\ln D}{p-1},
\end{align*} where we used the identity $$\sum_{k=0}^{j-1} (j-k)p^k=\frac{1}{p-1}\left(\frac{p^{j+1}-p}{p-1}-j\right).$$ Let $j_0=j_0(n,\ell,\mu,\nu^2,p)\in\mathbb{N}$ be the smallest integer such that $$j_0\geqslant \frac{\ln D}{2\ln p}-\frac{p}{p-1}.$$  Then, for any $j\geqslant j_0$ it results
\begin{align}\label{lower bound log Cj}
\ln C_j & \geqslant p^j\left(\ln C_0-\frac{2p\ln p}{(p-1)^2}+\frac{\ln D}{p-1}\right)  = p^j \ln \big( \widetilde{D}\varepsilon^p\big),
\end{align} where $\widetilde{D}\doteq K D^{1/(p-1)}p^{-(2p)/(p-1)^2}$.

Finally, from \eqref{lowe bound U j} we can prove that $U$ blows up in finite time since the lower bound diverges as $j\to \infty$ and, besides, we can derive an upper bound estimate for the lifespan of the solution. Combining \eqref{alpha j}, \eqref{beta j} and \eqref{lower bound log Cj}, for $j\geqslant j_0$ and $t\geqslant T_1$ we obtain
\begin{align*}
U(t) & \geqslant \exp\left(p^j \ln \big( \widetilde{D}\varepsilon^p\big)\right) t^{-\alpha_j} (t-T_1)^{\beta_j}  \\
& \geqslant \exp\left(p^j \left( \ln \big( \widetilde{D}\varepsilon^p\big)-\Big(\tfrac{r_2+1}{p-1}+n(\ell+1)+\alpha_0\Big)\ln t +\Big(\tfrac{r_2+3}{p-1}+\beta_0\Big)\ln (t-T_1)\right)\right)  t^{\frac{r_2+1}{p-1}+n(\ell+1)} (t-T_1)^{-\frac{r_3+1}{p-1}} .
\end{align*} For $t>2T_1$ it holds the relation $\ln(t-T_1)\geqslant \ln t-\ln 2$, consequently, for $j\geqslant j_0$ we get
\begin{align*}
U(t) & \geqslant \exp\left(p^j \left( \ln \big( \widehat{D}\varepsilon^p\big)+\Big(\tfrac{2}{p-1}+\beta_0-\alpha_0-n(\ell+1)\Big)\ln t \right)\right) t^{\frac{r_2+1}{p-1}+n(\ell+1)} (t-T_1)^{-\frac{r_3+1}{p-1}} ,
\end{align*} where $\widehat{D}\doteq 2^{-\frac{r_2+3}{p-1}-\beta_0}\widetilde{D}$. 

Let us write explicitly the factor that multiplies $\ln t$ in the previous estimate
\begin{align*}
\tfrac{2}{p-1}+\beta_0-\alpha_0-n(\ell+1) &= \tfrac{2}{p-1} + (n-1)(\ell+1)+2- \Big[\tfrac{n-1}{2}(\ell+1)+\tfrac{\ell+\mu}{2}\Big]p -n(\ell+1) \\
& =  \tfrac{1}{p-1} \Big\{ -\Big[\tfrac{n-1}{2}(\ell+1)+\tfrac{\ell+\mu}{2}\Big]p(p-1)-(\ell+1)(p-1) +2p \Big\} \\
& =  \tfrac{1}{p-1} \Big\{ -\Big[\tfrac{n-1}{2}(\ell+1)+\tfrac{\ell+\mu}{2}\Big]p^2+\Big[\tfrac{n+1}{2}(\ell+1)+\tfrac{\mu-3\ell}{2}\Big]p+\ell+1\Big\} \\ &= \tfrac{\theta(n,\ell,\mu,p)}{p-1},
\end{align*}  where $\theta(n,\ell,\mu,p)$ is defined in \eqref{def theta} and is a positive quantity thanks to the condition $p<p_{\mathrm{Str}}(n+\frac{\mu}{\ell+1},\ell)$ on the exponent of the nonlinear term. Then, for $j\geqslant j_0$ and $t\geqslant 2T_1$ we arrived at
\begin{align}\label{final lower bound U}
U(t) & \geqslant \exp\left(p^j \left( \ln \left( \widehat{D}\varepsilon^p t^{\frac{\theta(n,\ell,\mu,p)}{p-1}} \right) \right)\right) t^{\frac{r_2+1}{p-1}+n(\ell+1)} (t-T_1)^{-\frac{r_3+1}{p-1}}.
\end{align}
Let us fix $\varepsilon_0=\varepsilon_0(u_0,u_1,n,p,\ell,\mu,\nu^2,R)$ such that $0<\varepsilon_0< (2T_1)^{-\theta(n,\ell,\mu,p)/(p(p-1))}\widehat{D}^{-1/p}$.  Then, for any $\varepsilon\in (0,\varepsilon_0]$ and any $t>\widehat{D}^{-(p-1)/\theta(n,\ell,\mu,p)}\varepsilon^{-p(p-1)/\theta(n,\ell,\mu,p)}$ we have
\begin{align*}
t>2T_1 \quad \mbox{and} \quad \ln \left( \widehat{D}\varepsilon^p t^{\frac{\theta(n,\ell,\mu,p)}{p-1}} \right) >0,
\end{align*} so, letting $j\to \infty$ in \eqref{final lower bound U}, the lower bound for $U$ blows up. Thus, we proved that $U$ is not finite for $t\gtrsim \varepsilon^{-p(p-1)/\theta(n,\ell,\mu,p)}$, that is, we showed the lifespan estimate in \eqref{lifespan est thm main} for $p<p_{\mathrm{Str}}(n+\frac{\mu}{\ell+1},\ell)$.

The case $1<p<p_{\mathrm{Fuj}}\big((\ell+1)n+\tfrac{\mu-1}{2}-\tfrac{\sqrt{\delta}}{2}\big)$ can be treated in a completely analogous way. Indeed, employing \eqref{lower bound U Fujita type} in place of \eqref{1st lower bound U Str}, so that, $C_0=I \varepsilon$ and $\beta_0-\alpha_0=\frac{1-\mu}{2}+\frac{\sqrt{\delta}}{2}$, we determine the following lower bound estimate for $U$ instead of \eqref{final lower bound U}
\begin{align*}
U(t) & \geqslant \exp\left(p^j \left( \ln \left( \widetilde{C}\varepsilon t^{\frac{2}{p-1}-(\ell+1)n+\beta_0-\alpha_0} \right) \right)\right) t^{\frac{r_2+1}{p-1}+n(\ell+1)} (t-T_1)^{-\frac{r_3+1}{p-1}}
\end{align*} for any $t\geqslant 2T_1$ and for $j$ greater than a suitable $j_1(n,\ell,\mu,\nu^2,p)\in\mathbb{N}$, where $\widetilde{C}$ is a certain positive constant. Thanks to the assumption on $p$, from this estimate we conclude the validity of the second upper bound estimate in \eqref{lifespan est thm main} by repeating the same argument as in the first case.

\begin{remark} In the case $\delta=0$ and for $p<p_{\mathrm{Fuj}}((\ell+1)n+r_1)$, from \eqref{U lin def} we see that we can actually improve the lower bound estimate \eqref{lower bound U Fujita type} by a logarithmic factor. By using a slicing procedure as in the next section, it is possible to prove the following upper bound estimate for the lifespan
\begin{align*}
T(\varepsilon)^{\frac{2}{p-1}-((\ell+1)n+r_1)} \ln T(\varepsilon)\lesssim \varepsilon^{-1}.
\end{align*}
\end{remark}

\section{Proof of Theorem \ref{Thm crit Fuj}}

In order to prove Theorem \ref{Thm crit Fuj} we derive a sequence of lower bound estimates for $U(t)$ with additional logarithmic factors. According to this purpose, we introduce  $\{\ell_j\}_{j\in\mathbb{N}}$ such that $\ell_{j}\doteq 2-2^{-(j+1)}$. We are going to use this sequence to apply a \emph{slicing procedure} in the iteration argument, following the ideas introduced in the paper \cite{AKT00}.  More precisely, we establish the following estimates
\begin{align}\label{U(t) l.b. crit Fuj}
U(t)\geqslant K_j t^{-r_1} \left(\ln\left(\frac{t}{\ell_j}\right)\right)^{\gamma_j}
\end{align} for any $j\in\mathbb{N}$ and any $t\geqslant \ell_j$, where $\{K_j\}_{j\in\mathbb{N}}$, $\{\gamma_j\}_{j\in\mathbb{N}}$ are sequences of nonnegative numbers to be determined iteratively. From \eqref{U lin def} we obtain that
\begin{align*}
U(t)\geqslant \begin{cases} I \,\varepsilon t^{-r_1} & \mbox{if} \ \delta>0, \\ I \,\varepsilon t^{-r_1} \ln t & \mbox{if} \ \delta=0, \end{cases}
\end{align*} which implies the validity of \eqref{U(t) l.b. crit Fuj} for $j=0$ provided that $K_0\doteq I\varepsilon$ and $$\gamma_0\doteq \begin{cases} 0 & \mbox{if} \ \delta>0, \\
1 & \mbox{if} \ \delta=0. \end{cases}$$
Let us proceed with the inductive step. We plug \eqref{U(t) l.b. crit Fuj} for some $j\in\mathbb{N}$ in \eqref{iteration frame for U} and we prove the validity of \eqref{U(t) l.b. crit Fuj} for $j+1$, prescribing suitable recursive relations for the terms $K_{j+1}$ and $\gamma_{j+1}$. Therefore, for $t\geqslant \ell_{j}$ we have
\begin{align*}
U(t) & \geqslant C K_j^p t^{-r_1} \int_{\ell_j}^t s^{r_1-r_2-1} \int_{\ell_j}^s \tau^{r_2+1-n(\ell+1)(p-1)-r_1 p} \left(\ln\left(\tfrac{\tau}{\ell_j}\right)\right)^{p \gamma_j}\, \mathrm{d}\tau \, \mathrm{d}s \\
&\geqslant C K_j^p t^{-r_1} \int_{\ell_j}^t s^{r_1-r_2-1 -[(\ell+1)n+r_1](p-1)} \int_{\ell_j}^s \tau^{r_2-r_1+1} \left(\ln\left(\tfrac{\tau}{\ell_j}\right)\right)^{p \gamma_j}\, \mathrm{d}\tau \, \mathrm{d}s. 
\end{align*} For $t\geqslant \ell_{j+1}$, we can shrink the interval of integration in the $\tau$-integral as follows
\begin{align*}
U(t) & \geqslant C K_j^p t^{-r_1} \int_{\ell_{j+1}}^t s^{r_1-r_2-1 -[(\ell+1)n+r_1](p-1)} \int_{\tfrac{\ell_j \, s}{\ell_{j+1}}}^s \tau^{r_2-r_1+1} \left(\ln\left(\tfrac{\tau}{\ell_j}\right)\right)^{p \gamma_j}\, \mathrm{d}\tau \, \mathrm{d}s \\
& \geqslant C K_j^p t^{-r_1} \int_{\ell_{j+1}}^t s^{r_1-r_2-1 -[(\ell+1)n+r_1](p-1)} \left(\ln\left(\tfrac{s}{\ell_{j+1}}\right)\right)^{p \gamma_j}\int_{\tfrac{\ell_j \, s}{\ell_{j+1}}}^s \left(\tau-\tfrac{\ell_j \, s}{\ell_{j+1}}\right)^{r_2-r_1+1} \, \mathrm{d}\tau \, \mathrm{d}s \\
& =  C (r_2-r_1+2)^{-1} \left(1-\tfrac{\ell_j}{\ell_{j+1}}\right)^{r_2-r_1+2} K_j^p t^{-r_1} \int_{\ell_{j+1}}^t s^{1 -[(\ell+1)n+r_1](p-1)} \left(\ln\left(\tfrac{s}{\ell_{j+1}}\right)\right)^{p \gamma_j} \, \mathrm{d}s. 
\end{align*} Due to $p=p_{\mathrm{Fuj}}\left((\ell+1)n+r_1\right)$, the power of $s$ is actually $-1$ in the last integral, so, for $t\geqslant \ell_{j+1}$ we obtain
\begin{align*}
U(t) & \geqslant C (r_2-r_1+2)^{-1} \left(1-\tfrac{\ell_j}{\ell_{j+1}}\right)^{r_2-r_1+2} K_j^p t^{-r_1} \int_{\ell_{j+1}}^t s^{-1} \left(\ln\left(\tfrac{s}{\ell_{j+1}}\right)\right)^{p \gamma_j} \, \mathrm{d}s \\
&=  C (r_2-r_1+2)^{-1} \left(1-\tfrac{\ell_j}{\ell_{j+1}}\right)^{r_2-r_1+2} K_j^p (p\gamma_j+1)^{-1} t^{-r_1} \left(\ln\left(\tfrac{t}{\ell_{j+1}}\right)\right)^{p \gamma_j+1}, 
\end{align*} which is exactly \eqref{U(t) l.b. crit Fuj} for $j+1$, provided that 
\begin{align*}
\gamma_{j+1} & =p\gamma_j+1, \\
K_{j+1} &=C (r_2-r_1+2)^{-1} \left(1-\tfrac{\ell_j}{\ell_{j+1}}\right)^{r_2-r_1+2}  (p\gamma_j+1)^{-1} K_j^p.
\end{align*}
Applying iteratively the recursive relation $\gamma_j=p\gamma_{j-1}+1$, we obtain
\begin{align}
\gamma_j & = p^j\gamma_0+\sum_{k=0}^{j-1}p^k=\left(\gamma_0+\tfrac{1}{p-1}\right)p^j-\tfrac{1}{p-1} \label{estimate gamma j}
\\ &\leqslant \left(\gamma_0+\tfrac{1}{p-1}\right)p^j. \notag
\end{align} Moreover, $1-\frac{\ell_{j-1}}{\ell_{j}}>2^{-(j+2)}$. Consequently, combining the previous considerations, for any $j\geqslant 1$ we have
\begin{align*}
K_j\geqslant M Q^{-j}K_{j-1}^p,
\end{align*} where $M\doteq C (r_2-r_1+2)^{-1}2^{-2(r_2-r_1+2)}\left(\gamma_0+\frac{1}{p-1}\right)^{-1}$ and $Q\doteq 2^{r_2-r_1+2} p$. Applying the logarithmic function to both sides of the previous inequality and employing recursively the obtained estimate, we get
\begin{align*}
\ln K_j & \geqslant p \ln K_{j-1}-j\ln Q+\ln M \\
& \geqslant p^2 \ln K_{j-2}-(j+(j-1)p)\ln Q+(1+p) \ln M  \\ & \geqslant \ldots \geqslant p^j \ln K_0 - \ln Q\left( \, \sum_{k=0}^{j-1} (j-k)p^k\right)+\ln M \sum_{k=0}^{j-1}p^k \\
&= p^j\left(\ln K_0-\frac{p\ln Q}{(p-1)^2}+\frac{\ln M}{p-1}\right)+\frac{\ln Q}{p-1}\, j+\frac{p\ln Q}{(p-1)^2}-\frac{\ln M}{p-1}.
\end{align*} Let $j_2=j_2(n,\mu,\nu^2,\ell)\in\mathbb{N}$ be the smallest integer such that $j_1\geqslant \frac{\ln M}{\ln Q}-\frac{p}{p-1}$. Therefore, for any $j\in\mathbb{N}$, $j\geqslant j_2$ it holds
\begin{align}\label{estimate Kj}
\ln K_j\geqslant p^j \ln(\widetilde{M} \varepsilon),
\end{align} where $\widetilde{M}\doteq I M^{1/(p-1)}Q^{-p/(p-1)^2}$.

Finally, we combine \eqref{U(t) l.b. crit Fuj}, \eqref{estimate gamma j} and \eqref{estimate Kj}, so for any $j\geqslant j_2$ and any $t\geqslant 2\geqslant \ell_j$ we find
\begin{align*}
U(t) & \geqslant  K_j t^{-r_1} \left(\ln\left(\tfrac{t}{2}\right)\right)^{\left(\gamma_0+\frac{1}{p-1}\right)p^j-\frac{1}{p-1}} \\
& = \exp\left[ \ln K_j+\left(\gamma_0+\tfrac{1}{p-1}\right)p^j\ln\left(\ln\left(\tfrac t2\right)\right) \right] t^{-r_1} \left(\ln\left(\tfrac{t}{2}\right)\right)^{-\frac{1}{p-1}} \\
& \geqslant \exp\left[ p^j \left(\ln (\widetilde{M}\varepsilon)+\left(\gamma_0+\tfrac{1}{p-1}\right)\ln\left(\ln\left(\tfrac t2\right)\right) \right)\right] t^{-r_1} \left(\ln\left(\tfrac{t}{2}\right)\right)^{-\frac{1}{p-1}} \\
& = \exp\left[ p^j\ln\left(\widetilde{M}\varepsilon\left(\ln\left(\tfrac t2\right)\right)^{\gamma_0+1/(p-1)}\right) \right] t^{-r_1} \left(\ln\left(\tfrac{t}{2}\right)\right)^{-\frac{1}{p-1}}.
\end{align*} Since $\ln (\tfrac {t}{2})\geqslant \tfrac{1}{2} \ln t$ for $ t\geqslant 4$, for $j\geqslant j_2$ and $t\geqslant 4$ we arrive at 
\begin{align}\label{final lb est U(t) crit Fuj}
U(t) & \geqslant  \exp\left[ p^j\ln\left(\widehat{M}\varepsilon\left(\ln t\right)^{\gamma_0+1/(p-1)}\right) \right] t^{-r_1} \left(\ln\left(\tfrac{t}{2}\right)\right)^{-\frac{1}{p-1}}, 
\end{align} where $\widehat{M}\doteq 2^{-(\gamma_0+1/(p-1))} \widetilde{M}$.

 Let us fix $\varepsilon_0=\varepsilon(u_0,u_1,n,\ell,\mu,\nu^2,R)>0$ such that $\varepsilon_0\leqslant (2\ln 2)^{-(\gamma_0+1/(p-1))}\widehat{M}^{-1}$. Then, for any $\varepsilon\in (0,\varepsilon_0]$ and any $t>\exp\left((\widehat{M}\varepsilon)^{-(\gamma_0+1/(p-1))^{-1}}\right)$ we have
\begin{align*}
t\geqslant 4  \quad \mbox{and} \quad \ln\left(\widehat{M}\varepsilon\left(\ln t\right)^{\gamma_0+1/(p-1)}\right)>0,
\end{align*} thus, letting $j\to \infty$ in \eqref{final lb est U(t) crit Fuj} we obtain that the lower bound for $U(t)$ blows up. Hence, we proved that for $\ln t \gtrsim \varepsilon^{-(\gamma_0+1/(p-1))^{-1}}$ the average $U(t)$ may not be finite. This completes the proof and shows the upper bound estimate \eqref{lifespan est thm crit Fuj} for the lifespan.

\section{Concluding remarks}

In this final section, we analyze the result obtained in Theorem \ref{Thm main}. Setting
\begin{align*}
p_{\mathrm{c}}(n,\ell,\mu,\nu^2)\doteq \max\left\{p_{\mathrm{Str}}\!\left(n+\tfrac{\mu}{\ell+1},\ell\right),p_{\mathrm{Fuj}}\!\left((\ell+1)n+\tfrac{\mu-1}{2}-\tfrac{\sqrt{(\mu-1)^2-4\nu^2}}{2}\right) \right\},
\end{align*} from Theorem \ref{Thm main} we know that a blow-up result holds for $1<p<p_{\mathrm{c}}(n,\ell,\mu,\nu^2)$ provided that $\delta \geqslant 0$ and that the compactly supported Cauchy data fulfill suitable sign assumptions. As in the case of the semilinear wave equation with scale-invariant damping and mass, the condition $\delta \geqslant 0$ implies that the damping term is dominant over the mass one. From a technical viewpoint, this assumption guarantees the possibility to use the double multiplier technique (cf. \cite{PT18,ITW21}) while deriving the iteration frame \eqref{iteration frame for U}.

Combining the blow-up result from this paper with the global existence result for small data solutions in \cite[Corollary 5.1]{DAbb21}, it follows that $p_{\mathrm{c}}(n,\ell,\mu,\nu^2)$ is the critical exponent for \eqref{semilinear EPDT} when $n=1$.

Moreover, it is reasonable to conjecture that the exponent $p_{\mathrm{c}}(n,\ell,\mu,\nu^2)$ is critical even for higher dimensions. Indeed, for $\mu=\nu^2=0$ we have that $p_{\mathrm{c}}(n,\ell,0,0)=p_{\mathrm{Str}}(n,\ell)$ for $n\geqslant 2$ and $p_{\mathrm{c}}(n,\ell,0,0)=p_{\mathrm{Fuj}}(\ell)$ for $n=1$ (see \cite[Remark 1.6]{HWY20} for the one-dimensional case) according to the results for \eqref{semilinear Tricomi} with power nonlinearity that we recalled in the introduction. In particular, when $\ell=0$ too we find that $p_{\mathrm{c}}(n,0,0,0)$ is the solution to the quadratic equation $(n-1)p^2-(n+1)p-2=0$, namely, the celebrated exponent named after the author of \cite{Str81} which is the critical exponent for the semilinear wave equation. 

Furthermore, we point out that for $\ell=0$ and for $\nu^2=0,\ell\in (-1,0)$ the exponent $p_{\mathrm{c}}(n,\ell,\mu,\nu^2)$ coincides with \eqref{critical exponent scale invariant} and \eqref{critical exponent EdS}, respectively. 

In Theorem \ref{Thm crit Fuj}, the blow-up of local solutions is proved (under suitable assumptions for the Cauchy data) in the critical case with Fujita-type exponent. In the forthcoming paper \cite{LPT24}, the blow-up in the other critical case of Strauss-type is considered, with the approach for the critical classical semilinear wave equation developed in \cite{TakWak11}.

\section*{Acknowledgments}

A. Palmieri is member of the \emph{Gruppo Nazionale per L’Analisi Matematica, la Probabilità e le loro Applicazioni} (GNAMPA) of the \emph{Instituto Nazionale di Alta Matematica} (INdAM). A. Palmieri has been partially supported by  INdAM - GNAMPA Project 2024 ``Modelli locali e non-locali con perturbazioni non-lineari'' CUP E53C23001670001 and by ERC Seeds UniBa Project ``NWEinNES'' CUP H93C23000730001.

\end{document}